\documentclass[twoside,a4paper,12pt,centertags,reqno]{amsart} 

\usepackage{amsmath,amssymb,verbatim,vmargin}
\usepackage{color}
\usepackage{tikz}
\usepackage{tikz-cd}
\usetikzlibrary{matrix,calc}
\usepackage{hyperref}
\hypersetup{colorlinks,bookmarks=true,linktocpage=true,citecolor=blue,linkcolor=magenta}

\usepackage{eulervm}   

\allowdisplaybreaks 
\usepackage{color}
\usepackage[all]{xy} 

\usepackage{mathtools}
\usepackage{MnSymbol} 

\newtheorem{thm}{Theorem}[section]

\newtheorem{ex}[thm]{Example}

\usepackage{enumerate}
\usepackage{enumitem}

\newcommand{\pd}{\partial}

\newcommand{\bC}{{\mathbb C}}

\newcommand{\bD}{{\mathbb D}}

\newcommand{\bZ}{{\mathbb Z}}

\newcommand{\cH}{{\mathcal H}}

\def\barint_#1{\mathchoice
            {\mathop{\vrule width 6pt
height 3 pt depth -2.5pt
                    \kern -9.5pt
\intop \kern -4pt}\nolimits_{#1}}%
            {\mathop{\vrule width 5pt height
3 pt depth -2.6pt
                    \kern -6.5pt
\intop \kern -4pt}\nolimits_{#1}}%
            {\mathop{\vrule width 5pt height
3 pt depth -2.6pt
                    \kern -6pt
\intop \kern -4pt}\nolimits_{#1}}%
            {\mathop{\vrule width 5pt height
3 pt depth -2.6pt
          \kern -6pt \intop \kern -4pt}\nolimits_{#1}}}
          
           \def\bariint_#1{\mathchoice
            {\mathop{\vrule width 15pt
height 3 pt depth -2.5pt
                    \kern -15.8pt
\intop \kern -8pt\intop \kern -4pt}\nolimits_{#1}}%
            {\mathop{\vrule width 9pt height
3 pt depth -2.6pt
                    \kern -10.5pt
\intop \kern -8pt\intop \kern -4pt}\nolimits_{#1}}%
            {\mathop{\vrule width 9pt height
3 pt depth -2.6pt
                    \kern -10pt
\intop \kern -8pt\intop \kern -4pt}\nolimits_{#1}}%
            {\mathop{\vrule width 9pt height
3 pt depth -2.6pt
          \kern -8pt \intop \kern -10pt\intop \kern -4pt}
      \nolimits_{  #1}}}

\def\barintlim_#1{\mathchoice
            {\mathop{\vrule width 6pt
height 3 pt depth -2.5pt
                    \kern -8.8pt
\intop \kern -4pt}\limits_{#1}}%
            {\mathop{\vrule width 5pt height
3 pt depth -2.6pt
                    \kern -6.5pt
\intop \kern -4pt}\limits_{#1}}%
            {\mathop{\vrule width 5pt height
3 pt depth -2.6pt
                    \kern -6pt
\intop \kern -4pt}\limits_{#1}}%
            {\mathop{\vrule width 5pt height
3 pt depth -2.6pt
          \kern -6pt \intop \kern -4pt}\limits_{#1}}}
          
           \def\bariintlim_#1{\mathchoice
            {\mathop{\vrule width 15pt
height 3 pt depth -2.5pt
                    \kern -15.8pt
\intop \kern -8pt\intop \kern -4pt}\limits_{#1}}%
            {\mathop{\vrule width 9pt height
3 pt depth -2.6pt
                    \kern -10.5pt
\intop \kern -8pt\intop \kern -4pt}\limits_{#1}}%
            {\mathop{\vrule width 9pt height
3 pt depth -2.6pt
                    \kern -10pt
\intop \kern -8pt\intop \kern -4pt}\limits_{#1}}%
            {\mathop{\vrule width 9pt height
3 pt depth -2.6pt
          \kern -8pt \intop \kern -10pt\intop \kern -4pt}
      \limits_{  #1}}}
          
\renewcommand{\iint}{\int \kern -8pt\int}       

\begin{document}


\title[Number-phase uncertainty via Bergman spaces]
{Minimum number-phase uncertainty states via weighted Bergman spaces}

\author{Yi C. Huang}
\address{School of Mathematical Sciences, Nanjing Normal University, Nanjing 210023, People's Republic of China}
\email{Yi.Huang.Analysis@gmail.com}
\urladdr{https://orcid.org/0000-0002-1297-7674}

\date{\today}
\begin{abstract}
The number-phase uncertainty result of Luo via the Hardy space on unit disc (Phys Lett A, 2000) is extended in this paper to the scale of weighted Bergman spaces.
The minimum uncertainty states are thereby explicitly identified.
\end{abstract}


\maketitle

\section{Introduction}

We are interested in the number-phase observable on the quantum mechanical Hilbert space $\cH$ of harmonic oscillator with one degree of freedom.
The (complex) inner product for $\cH$ is denoted by $\langle\cdot,\cdot\rangle$.
As usual, we shall adopt Dirac's bra-ket notations $|\phi\rangle$ and $\langle\psi|$
to represent the vector $\phi$ in $\cH$ and the linear functional $$\langle\psi|: \phi\mapsto \langle\psi|\phi\rangle$$ that acts on $\cH$.
Let $N$ be the number operator with normalised eigenstates denoted by $|n\rangle$, $n=0,1,2,\cdots$.
Then $N$ can be expressed formally as
$$N=\sum_{n=0}^\infty n|n\rangle\langle n|.$$
Let $\Phi$ be the exponential phase operator proposed by Dirac \cite{Dir27} 
(see also Susskind and Glogower \cite{SusGlo64}, L\'evy-Leblond \cite{LevLeb76} and Newton \cite{New80} and the references therein)
$$\Phi=\sum_{n=0}^\infty |n\rangle\langle n+1|.$$
Note that $\Phi$ annihilates $|0\rangle$, and sends $|n+1\rangle$ to $|n\rangle$ for $n=0,1,2,\cdots$.
Apparently, 
\begin{equation} \label{e:leib}
[\Phi, N]=\Phi.
\end{equation}
The formal adjoint of $\Phi$ is given by
$$\Phi^*=\sum_{n=0}^\infty |n+1\rangle\langle n|.$$
Accordingly, we have the operator identity $\langle 0|\Phi^*=0$ on $\cH$.

As a concrete representation via analytic functions, 
we take $\cH$ as the weighted Bergman space $\cH_\lambda$ (see e.g. Hedenmalm, Korenblum and Zhu \cite{HedKorZhu00}) defined by
$$\cH_\lambda:=\bigg\{f:\bD\rightarrow\bC, \,\text{holomorphic},\,\langle f,f\rangle=\frac{\lambda-1}{\pi}\iint_{\bD}f(z)\overline{f(z)}(1-z\bar z)^{\lambda-2}dzd\bar z<\infty\bigg\}.$$
Here $\lambda>1$ is a weight parameter.
An orthonormal basis of $\cH_\lambda$ is
$$\bigg\{e_n(z):=\sqrt{\frac{\Gamma(n+\lambda)}{n!\Gamma(\lambda)}}z^n:n\geq0\bigg\}.$$
For a harmonic oscillator model based on $\cH_\lambda$, see for example Luo \cite{Luo97}.
The degenerate case $\lambda=1$ corresponds to the Hardy space in complex analysis 
$$\cH_1:=\bigg\{f:\bD\rightarrow\bC, \,\text{holomorphic},\,\langle f,f\rangle=\lim_{r\rightarrow1}\int_0^{2\pi}f(re^{i\theta})\overline{f(re^{i\theta})}\frac{d\theta}{2\pi}<\infty\bigg\},$$
with orthonormal (Taylor) basis $\{z^n:n\geq0\}$.
For further function theoretic studies on $\cH_1$ and $\{\cH_\lambda\}_{\lambda>1}$, see Garnett's monograph \cite{Gar81} and the aforementioned \cite{HedKorZhu00}.

\begin{ex}
In above analytic representation, we have
$$Nf(z)=z\frac{\pd}{\pd z}f(z),$$
and we denote by $\Phi_\lambda$ and $\Phi_\lambda^*$ the corresponding exponential phase operators.
For $\lambda=1$ we encounter with the backward/forward shifts (see Nikol'ski\u{\i}'s treatise \cite{Nik86})
$$\Phi_1 f(z)=\frac{f(z)-f(0)}{z},$$
$$\Phi_1^*f(z)=zf(z).$$
Hence, $\Phi_1 N=N$ and the Leibniz rule $$N\Phi_1=N-\Phi_1$$
leads to the commutation relation \eqref{e:leib}.
\end{ex}

Let $f\in\cH$ be a state with unit norm.
For any operator $A$ (not necessarily Hermitian) on $\cH$,
the expectation of $A$ in the state $f$ is defined as
$$\langle A\rangle=\langle A\rangle_f:=\langle f,Af\rangle.$$
The variance of $A$ in $f$ is then defined as
$$(\Delta A)^2=(\Delta_f A)^2:=\bigg\langle(A-\langle A\rangle)(A-\langle A\rangle)^*\bigg\rangle.$$
The minimum uncertainty states are the coherent states that minimise the uncertainty relation under investigation.
In this paper, we are interested in the number-phase uncertainty relation and we aim to minimise the quantity $(\Delta_f N)^2(\Delta_f \Phi_\lambda)^2$.

\begin{thm} \label{thm:main}
Let $\lambda>1$. 
The minimum uncertainty states for the number-phase pair $(N,\Phi)$ in $\cH_\lambda(\bD)$ can be parametrised as
\begin{equation} \label{e:para}
\begin{aligned}
(w,k)\in\bC\times\bZ_+  \mapsto f_{w,k}(z)=cz^k\sum_{n=0}^\infty \frac{w^n}{n!}e_n(z).
\end{aligned}
\end{equation}
Here $\bZ_+=\{0,1,2,\cdots\}$ and the normalisation constant $c\in\bC$ is determined by
\begin{equation} \label{e:norma}
I_{k,\lambda}(|w|^2)=(\bar cc)^{-1},
\end{equation}
where $$I_{k,\lambda}(t):=\sum_{n=0}^\infty \frac{t^n}{(n!)^2}\frac{(n+k)!}{n!}\frac{\Gamma(n+\lambda)}{\Gamma(n+k+\lambda)}.$$ 
In particular,
\begin{equation} \label{e:consist}
\langle N\rangle_{f_{w,k}}-w\overline{\langle \Phi_\lambda\rangle_{f_{w,k}}}=k.
\end{equation}
\end{thm}

Related results can be found in Carruthers and Nieto \cite{CarNie68}, Lerner, Huang and Walters \cite{LerHuaWal70} and Luo \cite{Luo00}.
For a nice survey exploring the number-phase statistics via analytic functions, see Vourdas \cite{Vou06} (and also the papers \cite{BriVouMan96, VouBriMan96}).

\section{Proof of Theorem \ref{thm:main}}

\medskip\textit{Derivation of \eqref{e:para}}.---Using Cauchy-Schwarz and noting that $N$ is Hermitian,
$$\begin{aligned}
(\Delta N)^2(\Delta \Phi_\lambda)^2&=\bigg\langle(N-\langle N\rangle)^*f,(N-\langle N\rangle)^*f\bigg\rangle\\
&\qquad\qquad\times\bigg\langle(\Phi_\lambda-\langle \Phi_\lambda\rangle)^*f,(\Phi_\lambda-\langle \Phi_\lambda\rangle)^*f\bigg\rangle\\
&\geq\bigg|\bigg\langle(N-\langle N\rangle)f,\left(\Phi^*_\lambda-\overline{\langle \Phi_\lambda\rangle}\right)f\bigg\rangle\bigg|^2.
\end{aligned}$$
In using Cauchy-Schwarz, the equality holds iff there exists $w\in\bC$ such that
$$(N-\langle N\rangle)f=w \left(\Phi_\lambda^*-\overline{\langle \Phi_\lambda\rangle}\right)f,$$
or by introducing $k=\langle N\rangle-w\overline{\langle \Phi_\lambda\rangle}$,
\begin{equation} \label{e:f}
Nf=w\Phi_\lambda^*f+kf.
\end{equation}
We can solve the operational part of \eqref{e:f}, 
$$Ng=w\Phi_\lambda^*g,$$ 
by the (normalised) eigenstate expansion, and
$$g(z)=g(0)\sum_{n=0}^\infty \frac{w^n}{n!}e_n(z).$$
Thus, we solve \eqref{e:f} with
$$f(z)=cz^k\sum_{n=0}^\infty \frac{w^n}{n!}e_n(z),$$
where $c$ is a normalisation constant, 
and for $f\in\cH_\lambda$ it is necessary that $k\in\bZ_+$.

\medskip\textit{Derivation of \eqref{e:norma}}.---Recall that $f$ has unit norm. Using $\|e_{n+k}\|_{\cH_\lambda}=1$ for all $n\geq0$,
$$\langle f,f\rangle=\bar cc\sum_{n=0}^\infty \frac{(w \overline w)^n}{(n!)^2}\frac{(n+k)!}{n!}\frac{\Gamma(n+\lambda)}{\Gamma(n+k+\lambda)}=1.$$
This gives \eqref{e:norma}. For convenience, let
$$G(n,k):=\frac{(n+k)!}{n!}\frac{\Gamma(n+\lambda)}{\Gamma(n+k+\lambda)}.$$
Note that for $\lambda=1$, $G(n,k)\equiv1$.
  
\medskip\textit{Derivation of \eqref{e:consist}}.---We compute
$$\begin{aligned}
\langle N\rangle_{f_{w,k}}
&=\bar cc\sum_{n=0}^\infty \frac{(w \overline w)^n}{(n!)^2}G(n,k)(n+k)\\
&=k+\bar cc\sum_{n=0}^\infty \frac{(w \overline w)^{n+1}}{(n!)(n+1)!}G(n+1,k).
\end{aligned}$$
Since $w\overline{\langle \Phi_\lambda\rangle}=\overline{\langle \overline w\Phi_\lambda\rangle}$ and
$$\begin{aligned}
\overline w\Phi_\lambda f(z)&=c (\overline w w)\sum_{n=0}^\infty \frac{w^n}{(n+1)!}\sqrt{G(n+1,k)}e_{n+k}(z),
\end{aligned}$$
we compute
$$w\overline{\langle \Phi_\lambda\rangle}=\bar cc\sum_{n=0}^\infty \frac{(w \overline w)^{n+1}}{n!(n+1)!}G(n+1,k).$$
Thus, the consistency equation \eqref{e:consist} is verified.

\section{Conclusion}

We extended S. Luo's number-phase uncertainty result \cite{Luo00} to the scale of weighted Bergman spaces.
His Hardy space result is the $\lambda\rightarrow1$ limit of Theorem \ref{thm:main}:
$$f_{w,k}(z) \rightarrow cz^k\sum_{n=0}^\infty \frac{w^n}{n!}z^n =cz^ke^{wz}.$$
which are shifted Barut-Girardello states, see e.g. Brif \cite{Bri95}.
The weighted Bergman spaces are useful in harmonic analysis, functional inequalities and quantum mechanical studies, 
see e.g. Luo \cite[Sect. 4-5]{Luo97} and Frank \cite{Fra23}. 
It would be (mathematically) interesting to consider the representation via other classes of analytic functions.

\bigskip

\section*{\textbf{Declarations}}

\bigskip

\textbf{Ethical Approval}. Not applicable.

\bigskip

\textbf{Declaration of competing interest}.
The author declares that he has no competing financial interests or personal relationships 
that could have appeared to influence the work reported in this paper.

\bigskip

\textbf{Availability of data and materials}.
Data sharing not applicable to this article as no data sets were generated or analysed during the current study.

\bigskip

\textbf{Funding}.
Research supported by National NSF grant of China (no. 11801274).

\bigskip

\textbf{Acknowledgement}.
The author would like to thank Professors Min LI (HUST), Yuan SUN (NJNU) and Haiyan XU (GDUT) for helpful communications.


\begin{thebibliography}{19}

\bibitem{Bri95} C. Brif,
Photon states associated with the Holstein-Primakoff realization of the $SU(1,1)$ Lie algebra,
\textit{Quantum Semiclas. Opt.} {\bf 7} (1995) 803--834.

\bibitem{BriVouMan96} C. Brif, A. Vourdas, A. Mann,
Analytic representations based on $SU(1,1)$ coherent states and their applications,
\textit{J. Phys. A} {\bf 29} (1996) 5873--5885.

\bibitem{CarNie68} P. Carruthers, M. M. Nieto,
Phase and angle variables in quantum mechanics,
\textit{Rev. Mod. Phys.} {\bf 40} (1968) 411--440.

\bibitem{Dir27} P. A. M. Dirac,
The quantum theory of the emission and absorption of radiation,
\textit{Proc. Roy. Soc. (London) A} {\bf 114} (1927) 243--265.

\bibitem{Fra23} R. L. Frank,
Sharp inequalities for coherent states and their optimizers, 
\textit{Adv. Nonlinear Stud.} {\bf 23} (2023) 20220050. (Special Issue: In honor of David Jerison)

\bibitem{Gar81} J. B. Garnett,
\textit{Bounded Analytic Functions}, Academic Press, New York, 1981.

\bibitem{HedKorZhu00} H. Hedenmalm. B. Korenblum, K. Zhu
\textit{Theory of Bergman Spaces}, Graduate Texts in Mathematics (GTM, volume {\bf 199}), 2000.

\bibitem{LerHuaWal70} E. C. Lerner, H. W. Huang, G. E. Walters,
Some mathematical properties of oscillator phase operators,
\textit{J. Math. Phys.} {\bf 11} (1970) 1679--1684.

\bibitem{LevLeb76} J. M. L\'evy-Leblond,
Who is afraid of nonhermitian operators? A quantum description of angle and phase,
\textit{Ann. Phys. (NY)} {\bf 101} (1976) 319--341.

\bibitem{Luo97} S. Luo,
A harmonic oscillator on the Poincar\'e disc and hypercontractivity,
\textit{J. Phys. A} {\bf 30} (1997) 5133--5139.

\bibitem{Luo00} S. Luo,
Minimum uncertainty states for Dirac's number-phase pair,
\textit{Phys. Lett. A} {\bf 275} (2000) 165--168.

\bibitem{New80} R. G. Newton,
Quantum action-angle variables for harmonic oscillators,
\textit{Ann. Phys.} {\bf 124} (1980) 327--346.

\bibitem{Nik86} N. K. Nikol'ski\u{\i}
\textit{Treatise on the Shift Operator (Spectral Function Theory)}, Grundlehren der Mathematischen Wissenschaften (GL, volume {\bf 273}), 1986.

\bibitem{SusGlo64} L. Susskind, J. Glogower,
Quantum mechanical phase and time operator,
\textit{Physics} {\bf 1} (1964) 49--61.

\bibitem{Vou06} A. Vourdas,
Analytic representations in quantum mechanics,
\textit{J. Phys. A} {\bf 39} (2006) R65--R141.

\bibitem{VouBriMan96} A. Vourdas, C. Brif, A. Mann,
Factorization of analytic representations in the unit disc and number-phase statistics of a quantum harmonic oscillator,
\textit{J. Phys. A} {\bf 29} (1996) 5887--5898.

\end{thebibliography}
\end{document}